УДК 519.856

# Универсальный метод поиска равновесий и стохастических равновесий в транспортных сетях

*Баймурзина Диляра Римовна*[1,2] dilyara.rimovna@gmail.com

*д.ф.-м.н. Гасников Александр Владимирович*[1,3] gasnikov.av@mipt.ru

*к.ф.-м.н. Гасникова Евгения Владимировна*[1] egasnikova@yandex.ru

*к.ф.-м.н. Двуреченский Павел Евгеньевич*[3,4] dvurechensky@iitp.ru

*Ершов Егор Иванович*[3] e.i.ershov@gmail.com

*Кубентаева Меруза Болатбековна*[1] kubikmeruza@yandex.ru

*Лагуновская Анастасия Александровна*[1] a.lagunovskaya@phystech.edu

[1] Факультет управления и прикладной математики Национального исследовательского Университета «Московский физико-технический институт».
141700, Россия, Московская область, г. Долгопрудный, Институтский переулок, д. 9
[2] Сколковский институт науки и технологий
143026, Россия, г. Москва, Территория Инновационного Центра “Сколково”, улица Нобеля, д. 3
[3] Институт проблем передачи информации им. А.А. Харкевича Российской академии наук
127051, Россия, г. Москва, Большой Каретный переулок, д.19 стр. 1
[4] Weierstrass Institute for Applied Analysis and Stochastics
10117, Germany, Berlin, Mohrenstr, 39

**Аннотация**

В работе предложен универсальный способ поиска обычных и стохастических равновесий в популяционных играх загрузки. Мы рассматриваем модели равновесного распределения потоков по путям Бэкмана и стабильной динамики. Поиск (стохастических) равновесий Нэша(–Вардропа) приводит к решению (энтропийно регуляризованных) задач выпуклой оптимизации. Данная работа посвящена тому, как эффективно решать такого рода задачи, а точнее двойственные к ним, с помощью недавно предложенного прямо-двойственного универсального градиентного метода, оптимально и адаптивно настраивающегося на гладкость решаемой задачи.

## 1. Введение

В работе [1] было анонсировано, что в цикле последующих статей планируется описать новые эффективные численные методы поиска равновесий (и стохастических равновесий [2]) в популярных моделях распределения транспортных потоков по путям (ребрам) графа транспортной сети [2–5]. В частности, это планировалось сделать для модели Бэкмана [3, 5] (называемой также BMW-моделью) и модели стабильной динамики [4, 5] (называемой также моделью Нестерова–деПальмы). В цикле работ [6–14] этот план начал осуществляться. Однако на данный момент остался открытым главный вопрос: какие из предложенных методов являются наилучшими? Одним из основных результатов данной работы является попытка сопоставления и унификации имеющихся здесь результатов с целью получения ответа на поставленный вопрос. Оказалось, что для поиска обычного (не стохастического) равновесия в модели Бэкмана наилучшим методом (по имеющимся на данный момент представлениям) является классический метод условного градиента Франк–Вульфа [3, 6], а для всех остальных случаев[1] разумно строить двойственную задачу[2] и решать ее универсальным градиентным методом Ю.Е. Нестерова [15] (см. также более поздние версии [16, 17]). Затем, пользуясь прямо-двойственностью этого метода [17], можно восстанавливать решение исходной (прямой задачи). Интересно отметить, что численные эксперименты (см. п. 7), проведенные с универсальным градиентным методом для поиска обычных (не стохастических) равновесий показали, что время работы метода пропорционально $\sim \tilde{\varepsilon}^{-1}$, где $\tilde{\varepsilon}$ – желаемая относительная точность (по функции) решения задачи, а не $\sim \tilde{\varepsilon}^{-2}$, как это можно было ожидать, исходя из теории решения негладких задач выпуклой оптимизации [18]. В отличие от специальных примеров, собранных в [15, 16], на которых наблюдался аналогичный эффект, в данной статье обнаружен, по-видимому, первый совершенно реальный такой пример с существенно негладкой функцией, имеющей огромное число всевозможных изломов.

В п. 2 приводятся экстремальные принципы, описывающие равновесную конфигурацию в транспортных сетях в четырех рассматриваемых случаях. Эти принципы формулируются как сепарабельные задачи выпуклой оптимизации (с композитом вида энтропии в случае поиска стохастических равновесий) при аффинных ограничениях. В основном в изложении используются результаты работ [1, 5, 8, 9, 11, 19–26].

В п. 3 с помощью достаточно стандартной техники выпуклого анализа строятся двойственные задачи, чтобы не проектироваться на аффинные ограничения при решении прямых задач (последняя операция является дорогой [18]). Во многом здесь используются результаты работ [8, 9, 12].

В п. 4, следуя [17], излагается специальный вариант универсального градиентного метода – универсальный метод подобных треугольников (УМПТ), в котором (в отличие от [15, 16]) используется только одна проекция на каждой итерации. Особое внимание уделяется прямо-двойственности этого метода. Это нашло отражение в необычной форме записи того, как сходится метод.

В п. 5 УМПТ применяется для решения двойственных задач, описанных в п. 3. Приводятся формулы восстановления решения прямой задачи. В основном изложение следует [9, 12]. Формулы восстановления решения задачи поиска равновесия в модели стабильной динамики, исходя из последовательности, сгенерированной УМПТ при решении двойственной задачи, по-видимому, приводятся впервые.

[1] Трех: поиск стохастических равновесий в модели Бэкмана и поиск обычных и стохастических равновесий в модели стабильной динамики.

[2] Эта задача оказывается негладкой при поиске обычных (не стохастических) равновесий.

В п. 5 внимание было сконцентрировано на том, как восстанавливать решения прямых задач, т.е. задач п. 2, исходя из приближенного решения двойственных задач, т.е. задач из п. 3. В п. 6 результаты п. 5 пополняются исследованием того, сколько итераций будет делать УМПТ (оценка констант метода) для достижения желаемой точности решения (прямой и двойственной задачи одновременно), и какая при этом будет стоимость одной итерации. Последняя задача сводится (при поиске стохастических равновесий) к вычислению значений и градиентов характеристических функций на транспортном графе [5, 9, 27] (для вычисления градиентов используется теория быстрого автоматического дифференцирования [28]) или к поиску кратчайших путей в транспортном графе [29] (при поиске обычных равновесий).

В п. 7 результаты всех предыдущих пунктов собираются вместе. Формулируется теорема, в которой описаны полученные в данной работе (верхние) теоретические оценки скорости сходимости УМПТ во всех четырех случаях. Приводятся результаты численных экспериментов. Эти результаты сопоставляются с ранее известными.

## 2. Экстремальные принципы для поиска равновесий в транспортных сетях

Рассмотрим транспортную сеть, которую будем представлять ориентированным графом $\langle V, E\rangle$, где $V$ – множество вершин (как правило, можно считать, что $|E|/4 \le |V| \le |E|$), а $E$ – множество ребер, $|E| = n$. Обозначим множество пар $w = (i, j)$ источник-сток через $OD$, $d_w$ – корреспонденция, отвечающая паре $w$, $x_p$ – поток по пути $p$; $P_w$ – множество путей, отвечающих корреспонденции $w$ (начинающихся в $i$ и заканчивающихся в $j$), $P = \bigcup_{w \in OD} P_w$ – множество всех путей. Затраты на прохождения ребра $e \in E$ описываются функцией $\tau_e(f_e)$, где $f_e$ – поток по ребру $e$.

Опишем марковскую логит-динамику (также говорят гиббсовскую динамику) в повторяющейся игре загрузки графа транспортной сети [19]. Пусть каждой корреспонденции отвечает $d_w M$ агентов ($M \gg 1$), $\tau_e(f_e) := \tau_e(f_e/M)$. Пусть имеется $TN$ шагов ($N \gg 1$). Каждый агент независимо от остальных на шаге $t+1$ выбирает с вероятностью

$$\frac{\lambda}{N} \frac{\exp\left(-G_p^t/\gamma\right)}{\sum_{q \in P_w} \exp\left(-G_q^t/\gamma\right)}.$$

путь $p \in P_w$, где $G_p^t$ – затраты на пути $p$ на шаге $t$ ($G_p^0 \equiv 0$), а с вероятность $1 - \lambda/N$ путь, который использовал на шаге $t$. Такая динамика отражает ограниченную рациональность агентов (см. п. 3, и часто используется в популяционной теории игр [19] и теории дискретного выбора [20]. Оказывается эта марковская динамика в пределе $N \to \infty$ превращается в марковскую динамику в непрерывном времени (вырождающуюся при $\gamma \to 0+$ в динамику наилучших ответов [19]), которая в свою очередь допускает два предельных перехода (обоснование перестановочности этих пределов см. в [21]): $T \to \infty$, $M \to \infty$ или $M \to \infty$, $T \to \infty$. При первом порядке переходов мы сначала ($T \to \infty$) согласно эргодической теореме для марковских процессов (в нашем случае марковский процесс – модель стохастической химической кинетики с унарными реакциями в условиях детального баланса [22–24]) приходим к финальной (=стационарной) вероятностной мере, имеющей в основе мультиномиальное распределение. С ростом числа агентов ($M \to \infty$) эта мера концентрируется около наиболее вероятного состояния, поиск которого сводится к решению задачи (1) ниже. Функционал в этой задаче оптимизации с точностью до потенцирования и мультипликативных и аддитивных констант соответствует исследуемой стационарной мере – то есть это функционал Санова [24, 25]. При обратном порядке переходов, мы сначала осуществляем, так называемый, канонический скейлинг [21, 22], приводящий к детерминированной кинетической динамике, описываемой СОДУ на $x$, а затем ($T \to \infty$) ищем аттрактор получившейся СОДУ. Глобальным аттрактором оказывается неподвижная точка, которая определяется решением задачи (1) ниже. Более того, функционал, стоящий в (1), является функцией Ляпунова полученной кинетической динамики (то есть функционалом Больцмана). Последнее утверждение – достаточно общий факт (функционал Санова, является функционалом Больцмана), верный при намного более общих условиях [24].

Итак, рассматривается следующая задача поиска стохастического равновесия Нэша–Вардропа в модели Бэкмана равновесного распределения транспортных потоков по путям (см., например, [1, 2]):

$$\sum_{e\in E}\sigma_e\left(f_e\right)+\gamma\sum_{w\in OD}\sum_{p\in P_w}x_p\ln\left(x_p/d_w\right)\to\min_{f=\Theta x,\,x\in X}, \tag{1}$$

где $\gamma>0$; $\sigma_e\left(f_e\right)=\int\limits_0^{f_e}\tau_e\left(z\right)dz$ – выпуклые функции;

$$\Theta=\left\|\delta_{ep}\right\|_{e\in E,\,p\in P}=\left\|\Theta^{\langle p\rangle}\right\|_{p\in P},\ \delta_{ep}=\begin{cases}1, & e\in p\\ 0, & e\notin p\end{cases};$$

$$X=\left\{x\geq 0:\ \sum_{p\in P_w}x_p=d_w,\ w\in OD\right\} \text{ – прямое произведение симплексов;}$$

в качестве $\tau_e\left(f_e\right)$ обычно выбирают BPR-функции [1, 3, 5, 30]:

$$\tau_e\left(f_e\right)=\bar{t}_e\cdot\left(1+\rho\cdot\left(f_e/\bar{f}_e\right)^4\right).$$

В пределе модели стабильной динамики [1]

$$\tau_e^{\mu}\left(f_e\right)=\bar{t}_e\cdot\left(1+\rho\cdot\left(f_e/\bar{f}_e\right)^{1/\mu}\right)\xrightarrow[\mu\to0+]{}\begin{cases}\bar{t}_e, & 0\leq f_e<\bar{f}_e\\ \left[\bar{t}_e,\infty\right), & f_e=\bar{f}_e\end{cases},$$

$$d\tau_e^{\mu}\left(f_e\right)/df_e\xrightarrow[\mu\to0+]{}0,\ 0\leq f_e<\bar{f}_e,$$

задача перепишется как [1]

$$\sum_{e\in E}f_e\bar{t}_e+\gamma\sum_{w\in OD}\sum_{p\in P_w}x_p\ln\left(x_p/d_w\right)\to\min_{\substack{f=\Theta x,\,x\in X\\ f\leq\bar{f}}}. \tag{2}$$

## 3. Переход к двойственной задаче

Запишем двойственную задачу к (1) [1, 6, 9] (далее мы используем обозначение $\operatorname{dom}\sigma^*$ – область определения сопряженной к $\sigma$ функции)

$$\min_{f,x}\left\{\sum_{e\in E}\sigma_e\left(f_e\right)+\gamma\sum_{w\in OD}\sum_{p\in P_w}x_p\ln\left(x_p/d_w\right):\ f=\Theta x,\ x\in X\right\}=$$

$$= \min_{f,x} \left\{ \sum_{e\in E} \max_{t_e \in \mathrm{dom}\,\sigma_e^*} \left[ f_e t_e - \sigma_e^*(t_e) \right] + \gamma \sum_{w\in OD} \sum_{p\in P_w} x_p \ln\left(x_p/d_w\right) :\ f = \Theta x,\ x\in X \right\} =$$

$$= \max_{t\in \mathrm{dom}\,\sigma^*} \left\{ \min_{f,x} \left[ \sum_{e\in E} f_e t_e + \gamma \sum_{w\in OD} \sum_{p\in P_w} x_p \ln\left(x_p/d_w\right) :\ f = \Theta x,\ x\in X \right] - \sum_{e\in E} \sigma_e^*(t_e) \right\} =$$

$$= -\min_{t\in \mathrm{dom}\,\sigma^*} \left\{ \gamma\psi\left(t/\gamma\right) + \sum_{e\in E} \sigma_e^*(t_e) \right\}, \qquad (3)$$

где[3]

$$\psi(t) = \sum_{w\in OD} d_w \psi_w(t),\ \psi_w(t) = \ln\left( \sum_{p\in P_w} \exp\left( -\sum_{e\in E} \delta_{ep} t_e \right) \right),$$

$$f = -\nabla\gamma\psi\left(t/\gamma\right),\ x_p = d_w \frac{\exp\left( -\frac{1}{\gamma} \sum_{e\in E} \delta_{ep} t_e \right)}{\sum_{q\in P_w} \exp\left( -\frac{1}{\gamma} \sum_{e\in E} \delta_{eq} t_e \right)},\ p \in P_w, \qquad (4)$$

для

$$\tau_e(f_e) = \bar{t}_e \cdot \left( 1 + \rho \cdot \left( \frac{f_e}{\bar{f}_e} \right)^{\frac{1}{\mu}} \right),$$

имеем [1] (BPR-функция, см. п. 2, получается при $\mu = 1/4$)

$$\sigma_e^*(t_e) = \sup_{f_e \ge 0} \left( (t_e - \bar{t}_e) \cdot f_e - \bar{t}_e \cdot \frac{\mu}{1+\mu} \cdot \rho \cdot \frac{f_e^{1+\frac{1}{\mu}}}{\bar{f}_e^{\frac{1}{\mu}}} \right) = \bar{f}_e \cdot \left( \frac{t_e - \bar{t}_e}{\bar{t}_e \cdot \rho} \right)^{\mu} \frac{(t_e - \bar{t}_e)}{1+\mu}.$$

Собственно, формула (3) есть не что иное, как отражение формулы $f = -\nabla\gamma\psi\left(t/\gamma\right)$ и связи $t_e = \tau_e(f_e)$, $e \in E$. Действительно, по формуле Демьянова–Данскина–Рубинова [31]

$$\frac{d\sigma_e^*(t_e)}{dt_e} = \frac{d}{dt_e} \max_{f_e \ge 0} \left\{ t_e f_e - \int_0^{f_e} \tau_e(z)\,dz \right\} = f_e :\ t_e = \tau_e(f_e).$$

В свою очередь, формула $f = -\nabla\gamma\psi\left(t/\gamma\right)$ может интерпретироваться, как следствие соотношений $f = \Theta x$ и формулы распределения Гиббса (логит-распределения)

[3] Обратим внимание, что вектор распределения потоков по путям $x$ при поиске стохастического равновесия ($\gamma > 0$) получается не разреженным в отличие от поиска обычного равновесия Нэша(–Вардропа) [6] ($\gamma \to 0+$). Как следствие, чтобы вычислить этот вектор требуются затраты существенно зависящие от потенциально огромной размерности $|P|$. К счастью, в приложениях, как правило, не требуется знание этого вектора, достаточно определить вектор потоков на ребрах $f$, который, как мы увидим ниже, может быть вычислен намного эффективнее (в частности, с затратами независящими от $|P|$).

$$x_p = d_w \frac{\exp\left(-\frac{1}{\gamma}\sum_{e\in E}\delta_{ep}t_e\right)}{\sum_{q\in P_w}\exp\left(-\frac{1}{\gamma}\sum_{e\in E}\delta_{eq}t_e\right)},\ p \in P_w.$$

При такой интерпретации связь задачи (3) с логит-динамикой, порождающей стохастические равновесия, наиболее наглядна. Последняя формула, в виду того, что $g_p(t) = \sum_{e\in E}\delta_{ep}t_e$ – затраты на пути $p$ на графе $\langle V,E\rangle$, ребра которого взвешены $t$, – есть ни что иное, как отражение следующего принципа поведения (ограниченной рациональности агентов [20]): каждый агент $k$ (пользователь транспортной сети), отвечающий корреспонденции $w \in OD$, выбирает маршрут следования $p \in P_w$, если

$$p = \arg\max_{q\in P_w}\left\{-g_q(t) + \xi_q^k\right\},$$

где независимые случайные величины $\xi_q^k$, имеют одинаковое двойное экспоненциальное распределение, также называемое распределением Гумбеля[4] [19, 20]:

$$P\left(\xi_q^k < \zeta\right) = \exp\left\{-e^{-\zeta/\gamma - E}\right\},\ \gamma > 0.$$

Отметим также, что если взять $E \approx 0.5772$ – константа Эйлера, то $\mathrm{M}\left[\xi_q^k\right] = 0$, $D\left[\xi_q^k\right] = \gamma^2\pi^2/6$. Распределение Гиббса получается в пределе, когда число агентов на каждой корреспонденции стремится к бесконечности (случайность исчезает и описание переходит на средние величины).

Сформулируем главный вывод из написанного выше. Если каждый пользователь сориентирован на вектор затрат $t$ на ребрах $E$, одинаковый для всех пользователей, и пытается выбрать кратчайший путь, исходя из зашумленной информации, то такое поведение пользователей в пределе, когда их число стремится к бесконечности $M \to \infty$, приводит к описанию распределения пользователей по путям/ребрам (4). *Равновесная конфигурация, описываемая решением задачи (3), характеризуется тем, что по вектору* $t$ *вычисляется согласно формуле (4) такой вектор* $f = \Theta x$, *что имеет место соотношение* $t = \left\{\tau_e\left(f_e\right)\right\}_{e\in E}$.

Заметим, что при $\gamma \to 0+$ распределение водителей по путям вырождается, и все водители (агенты) будут использовать только кратчайшие пути

$$-\lim_{\gamma\to 0+}\gamma\psi_w(t/\gamma) = \min_{p\in P_w} g_p(t), \tag{5}$$

$$-\lim_{\gamma\to 0+}\nabla\gamma\psi_w(t/\gamma) = \mathrm{conv}\left\{(0,1,\ldots,0,1)^T,\ldots,(1,1,\ldots,0,0)^T\right\}, \tag{6}$$

в правой части формулы (6) стоит субградиент (негладкой) функции длины кратчайшего пути, отвечающего корреспонденции $w$. Этот субградиент (в общем случае) есть множество, представимое в виде выпуклой оболочки векторов, отвечающих всевозможным кратчайшим путям, начинающихся в вершине $i$ и заканчивающихся в вершине $j$, где $w = (i,j)$, в рассматриваемом транспортном графе, ребра которого взвешены вектором $t$. Единицы стоят в векторах на местах, которые отвечают ребрам, входящим в данный кратчайший путь. Если кратчайший путь единственен, то субградиент превращается в обычный градиент. Далее (в УМПТ) можно выбирать произвольный элемент субградиента. От этого выбора приводимые ниже оценки не будут зависеть.

---

[4] Распределение Гумбеля можно объяснить исходя из идемпотентного аналога центральной предельной теоремы (вместо суммы случайных величин – максимум) для независимых случайных величин с экспоненциальным и более быстро убывающим правым хвостом [32]. Распределение Гумбеля возникает в данном контексте, например, если при принятии решения водитель собирает информацию с большого числа разных (независимых) зашумленных источников, ориентируясь на худшие прогнозы по каждому из путей.

Полезно также иметь в виду, что [19, 20]

$$\gamma\psi_w(t/\gamma) = \mathrm{M}_{\{\xi_p\}_{p\in P_w}}\left[\max_{p\in P_w}\{-g_p(t)+\xi_p\}\right].$$

В пределе стабильной динамики $\mu \to 0+$ задача (3) вырождается в задачу

$$\gamma\psi(t/\gamma) + \langle \bar{f}, t-\bar{t}\rangle \to \min_{t\ge\bar{t}}.$$

В пределе $\gamma \to 0+$ все приведенные выше формулы переходят в соответствующие формулы для классических (не стохастических) моделей Бэкмана и стабильной динамики, см., например, [6].

### 4. Универсальный метод подобных треугольников

Ниже, следуя [17], изложен возможный способ решения двойственной задачи (3), позволяющий при этом восстанавливать решение прямой задачи ((1) или (2)), исходя из подхода, описанного в [33, 34] (см., также первоисточники [35, 36]). Рассмотрим общую задачу выпуклой композитной оптимизации на множестве простой структуры

$$F(t) = \underbrace{\Phi(t)}_{\gamma\psi(t/\gamma)} + \underbrace{h(t)}_{\sum_{e\in E}\sigma_e^*(t_e)} \to \min_{t\in Q}. \tag{7}$$

Положим $R^2 = \frac{1}{2}\|t_* - y^0\|_2^2$, где $y^0 = \bar{t}$, $t_*$ – решение задачи (7) (если решение не единственно, то выбирается то, которое доставляет минимум $\|t_* - y^0\|_2^2$).

Пусть ($L_\nu \le \infty$, т.е. допускается равенство бесконечности)

$$\|\nabla\Phi(t) - \nabla\Phi(y)\|_2 \le L_\nu \|t-y\|_2^\nu,\ \nu\in[0,1],\ L_0<\infty. \tag{8}$$

Положим

$$\varphi_0(t) = \alpha_0\left[\Phi(y^0) + \langle\nabla\Phi(y^0), t-y^0\rangle + h(t)\right] + \frac{1}{2}\|t-y^0\|_2^2,$$

$$\varphi_{k+1}(t) = \varphi_k(t) + \alpha_{k+1}\left[\Phi(y^{k+1}) + \langle\nabla\Phi(y^{k+1}), t-y^{k+1}\rangle + h(t)\right],$$

Начнем описание универсального метода подобных треугольников (УМПТ) с самой первой итерации. Положим

$$A_0 = \alpha_0 = 1/L_0^0,\ k=0,\ j_0=0;\ t^0 = u^0 = \arg\min_{t\in Q}\varphi_0(t).$$

До тех пор пока

$$\Phi(t^0) > \Phi(y^0) + \langle\nabla\Phi(y^0), t^0-y^0\rangle + \frac{L_0^{j_0}}{2}\|t^0-y^0\|_2^2 + \frac{\alpha_0}{2A_0}\varepsilon,$$

выполнять

$$j_0 := j_0 + 1;\; L_0^{j_0} := 2^{j_0} L_0^0;\; t^0 := u^0 := \arg\min_{t\in Q} \varphi_0(t),\; (A_0 :=)\,\alpha_0 := \frac{1}{L_0^{j_0}}.$$

**<u>Универсальный Метод Подобных Треугольников</u>**

1. $L_{k+1}^0 = L_k^{j_k} \big/ 2$, $j_{k+1} = 0$.

2. $$\begin{cases} \alpha_{k+1} := \dfrac{1}{2L_{k+1}^{j_{k+1}}} + \sqrt{\dfrac{1}{4\left(L_{k+1}^{j_{k+1}}\right)^2} + \dfrac{A_k}{L_{k+1}^{j_{k+1}}}},\; A_{k+1} := A_k + \alpha_{k+1}; \\ y^{k+1} := \dfrac{\alpha_{k+1} u^k + A_k t^k}{A_{k+1}},\; u^{k+1} := \arg\min_{t\in Q} \varphi_{k+1}(t),\; t^{k+1} := \dfrac{\alpha_{k+1} u^{k+1} + A_k t^k}{A_{k+1}}. \end{cases} \quad (*)$$

До тех пор пока

$$\Phi\left(y^{k+1}\right) + \left\langle \nabla\Phi\left(y^{k+1}\right), t^{k+1} - y^{k+1} \right\rangle + \frac{L_{k+1}^{j_{k+1}}}{2} \left\| t^{k+1} - y^{k+1} \right\|_2^2 + \frac{\alpha_{k+1}}{2A_{k+1}} \varepsilon < \Phi\left(t^{k+1}\right),$$

выполнять

$$j_{k+1} := j_{k+1} + 1;\; L_{k+1}^{j_{k+1}} = 2^{j_{k+1}} L_{k+1}^0;\; (*).$$

3. $k := k+1$ и **go to** 1.

УМПТ сходится согласно оценке

$$A_N F\left(t^N\right) \le \min_{t\in Q} \left\{ \frac{1}{2} \left\| t - y^0 \right\|_2^2 + \sum_{k=0}^N \alpha_k \left[ \Phi\left(y^k\right) + \left\langle \nabla\Phi\left(y^k\right), t - y^k \right\rangle + h(t) \right] \right\}. \quad (9)$$

Отсюда

$$A_N F\left(t^N\right) \le \min_{t\in Q} \left\{ \frac{1}{2} \left\| t - y^0 \right\|_2^2 + \sum_{k=0}^N \alpha_k \left[ \Phi\left(y^k\right) + \left\langle \nabla\Phi\left(y^k\right), t - y^k \right\rangle + h(t) \right] \right\} + \frac{\varepsilon}{2} \le$$

$$\le \frac{1}{2} \left\| t_* - y^0 \right\|_2^2 + \sum_{k=0}^N \alpha_k \underbrace{\left[ \Phi\left(y^k\right) + \left\langle \nabla\Phi\left(y^k\right), t_* - y^k \right\rangle + h(t_*) \right]}_{\le \Phi(t_*) + h(t_*)} + \frac{\varepsilon}{2} \le$$

$$\le \frac{1}{2} \left\| t_* - y^0 \right\|_2^2 + \sum_{k=0}^N \alpha_k F(t_*) + \frac{\varepsilon}{2} = R^2 + A_N F(t_*) + \frac{\varepsilon}{2},$$

т.е.

$$F\left(t^N\right) - F(t_*) \le \frac{R^2}{A_N} + \frac{\varepsilon}{2}.$$

В [17] доказано, что $A_N \simeq 2R^2/\varepsilon$ при

$$N \simeq \inf_{\nu\in[0,1]}\left(\frac{L_\nu\cdot(16R)^{1+\nu}}{\varepsilon}\right)^{\frac{2}{1+3\nu}}, \tag{10}$$

т.е. при таком $N$ справедливо неравенство

$$F\left(t^N\right)-\min_{t\in Q}F\left(t\right)\le\varepsilon .$$

При этом среднее число вычислений значения функции на одной итерации будет $\approx 4$, а градиента функции $\approx 2$.

## 5. Приложение УМПТ к поиску равновесий в транспортных сетях

Задачи, получающиеся из (3) при $\mu>0$ и $\mu\to 0+$, можно не различать по сложности при композитном подходе к ним [17, 37]. В задаче, отвечающей $\mu\to 0+$, сепарабельный композит $h(t)$ проще сепарабельного композита задачи с $\mu>0$, но зато дополнительно добавляется сепарабельное ограничение простой структуры[5]. В любом случае, основные затраты на каждой итерации при использовании УМПТ в композитном варианте [17, 37] связаны с необходимостью расчета градиента гладкой функции $\gamma\psi(t/\gamma)$. Тому как это можно эффективно делать будет посвящен следующий пункт 6. Здесь же отметим, что константы Гёльдера (см. (8)) градиента этой функции в 2-норме могут быть оценены следующим образом [38– 40]

$$L_1\le\frac{1}{\gamma}\sum_{w\in OD}d_w\max_{p\in P_w}\left\|\Theta^{\langle p\rangle}\right\|_2^2\le\frac{Hd}{\gamma}, \tag{11}$$

$$L_0\le 2\sum_{w\in OD}d_w\max_{p\in P_w}\left\|\Theta^{\langle p\rangle}\right\|_2\le 2\sqrt{H}d \text{ (при } \gamma\to 0+\text{)}, \tag{12}$$

где $d=\sum_{w\in OD}d_w$, $H=\max_{p\in P}\left\|\Theta^{\langle p\rangle}\right\|_2^2$, т.е. $H$ – максимальное число ребер в пути. Как правило, можно считать, что $H=\mathrm{O}\left(\sqrt{n}\right)$ – диаметр графа с манхетенской структурой, то есть для квадратной решетки с $n$ ребрами.

В данном пункте внимание будет сосредоточено на том, как с помощью формулы (9) восстанавливать решение прямой задачи ((1) или (2)), исходя из последовательности, сгенерированной УМПТ при решении двойственной задачи (3). Описанные далее способы имеют много общего с методами, описанными в работах [12, 33, 34].

Положим

$$f^k=-\nabla\gamma\psi\left(y^k/\gamma\right),\ x_p^k=d_w\frac{\exp\left(-\frac{1}{\gamma}\sum_{e\in E}\delta_{ep}y_e^k\right)}{\sum_{q\in P_w}\exp\left(-\frac{1}{\gamma}\sum_{e\in E}\delta_{eq}y_e^k\right)},\ p\in P_w,\ w\in OD,$$

[5] Для BPR-функций ( $\mu=0.25$ ) все сводится к решению уравнения четвертой степени, см. п. 3 и [12].

$$\bar{f}^N = \frac{1}{A_N}\sum_{k=0}^{N} a_k f^k , \ \bar{x}^N = \frac{1}{A_N}\sum_{k=0}^{N} \alpha_k x^k .$$

Положим[6]

$$\left\{\tau_e\left(\bar{f}_e^N\right)\right\}_{e\in E} = \arg\min_{t\in \operatorname{dom}\sigma^*}\left\{\frac{1}{A_N}\left[\sum_{k=0}^{N}\alpha_k\cdot\left(\gamma\psi\left(y^k/\gamma\right)+\left\langle\nabla\gamma\psi\left(y^k/\gamma\right),t-y^k\right\rangle\right)\right]+\sum_{e\in E}\sigma_e^*\left(t_e\right)\right\}. \quad (13)$$

Из (14) следует, что

$$\gamma\psi\left(t^N/\gamma\right)+\sum_{e\in E}\sigma_e^*\left(t_e^N\right)\le$$

$$\le \min_{t\in \operatorname{dom}\sigma^*}\left\{\frac{1}{A_N}\left[\sum_{k=0}^{N}\alpha_k\cdot\left(\gamma\psi\left(y^k/\gamma\right)+\left\langle\nabla\gamma\psi\left(y^k/\gamma\right),t-y^k\right\rangle\right)\right]+\sum_{e\in E}\sigma_e^*\left(t_e\right)+\frac{1}{2A_N}\left\|t-\bar{t}\right\|_2^2\right\}+\frac{\varepsilon}{2}\le$$

$$\le \min_{t\in \operatorname{dom}\sigma^*}\left\{\frac{1}{A_N}\left[\sum_{k=0}^{N}\alpha_k\cdot\left(\gamma\psi\left(y^k/\gamma\right)+\left\langle\nabla\gamma\psi\left(y^k/\gamma\right),t-y^k\right\rangle\right)\right]+\sum_{e\in E}\sigma_e^*\left(t_e\right)\right\}+$$

$$+\frac{1}{A_N}\underbrace{\frac{1}{2}\sum_{e\in E}\left(\tau_e\left(\bar{f}_e^N\right)-\bar{t}_e\right)^2}_{\tilde{R}^2}+\frac{\varepsilon}{2},$$

следовательно

$$\gamma\psi\left(t^N/\gamma\right)+\sum_{e\in E}\sigma_e^*\left(t_e^N\right)-$$

$$-\min_{t\in \operatorname{dom}\sigma^*}\left\{\frac{1}{A_N}\left[\sum_{k=0}^{N}\alpha_k\cdot\left(\gamma\psi\left(y^k/\gamma\right)+\left\langle\nabla\gamma\psi\left(y^k/\gamma\right),t-y^k\right\rangle\right)\right]+\sum_{e\in E}\sigma_e^*\left(t_e\right)\right\}\le\frac{\tilde{R}^2}{A_N}+\frac{\varepsilon}{2}.$$

Учитывая, что

$$-\min_{t\in \operatorname{dom}\sigma^*}\left\{\frac{1}{A_N}\left[\sum_{k=0}^{N}\alpha_k\left\langle\nabla\gamma\psi\left(y^k/\gamma\right),t\right\rangle\right]+\sum_{e\in E}\sigma_e^*\left(t_e\right)\right\}=$$

$$=\max_{t\in \operatorname{dom}\sigma^*}\left\{\left\langle\frac{1}{A_N}\sum_{k=0}^{N}\alpha_k f^k,t\right\rangle-\sum_{e\in E}\sigma_e^*\left(t_e\right)\right\}=\sum_{e\in E}\sigma_e\left(\bar{f}_e^N\right),$$

$$-\frac{1}{A_N}\sum_{k=0}^{N}\alpha_k\cdot\left(\gamma\psi\left(y^k/\gamma\right)-\left\langle\nabla\gamma\psi\left(y^k/\gamma\right),y^k\right\rangle\right)=$$

---

[6] Минимум в (13) достигается во внутренней точки множества $\operatorname{dom}\sigma^*$, которое для BPR-функций имеет вид $t\ge\bar{t}$ , и лишь при $\mu\to 0+$ минимум может выходить на границу $t_e=\bar{t}_e$ .

$$= \frac{1}{A_N} \sum_{k=0}^{N} \alpha_k \cdot \left( \left\langle f^k, y^k \right\rangle + \gamma \sum_{w \in OD} \sum_{p \in P_w} x_p^k \ln\left( x_p^k / d_w \right) - \left\langle f^k, y^k \right\rangle \right) =$$

$$= \gamma \frac{1}{A_N} \sum_{k=0}^{N} \alpha_k \sum_{w \in OD} \sum_{p \in P_w} x_p^k \ln\left( x_p^k / d_w \right) \ge \gamma \sum_{w \in OD} \sum_{p \in P_w} \bar{x}_p^N \ln\left( \bar{x}_p^N / d_w \right),$$

получаем следующую оценку сверху на зазор двойственности

$$0 \le \left\{ \gamma\psi\left( t^N/\gamma \right) + \sum_{e \in E} \sigma_e^*\left( t_e^N \right) \right\} + \left\{ \sum_{e \in E} \sigma_e\left( \bar{f}_e^N \right) + \gamma \sum_{w \in OD} \sum_{p \in P_w} \bar{x}_p^N \ln\left( \bar{x}_p^N / d_w \right) \right\} \le$$

$$\le \left\{ \gamma\psi\left( t^N/\gamma \right) + \sum_{e \in E} \sigma_e^*\left( t_e^N \right) \right\} + \left\{ \sum_{e \in E} \sigma_e\left( \bar{f}_e^N \right) + \gamma \frac{1}{A_N} \sum_{k=0}^{N} \alpha_k \sum_{w \in OD} \sum_{p \in P_w} x_p^k \ln\left( x_p^k / d_w \right) \right\} \le \frac{\tilde{R}^2}{A_N} + \frac{\varepsilon}{2}, \quad (14)$$

следовательно

$$\sum_{e \in E} \sigma_e\left( \bar{f}_e^N \right) + \gamma \sum_{w \in OD} \sum_{p \in P_w} \bar{x}_p^N \ln\left( \bar{x}_p^N / d_w \right) - \left( \sum_{e \in E} \sigma_e\left( f_e^* \right) + \gamma \sum_{w \in OD} \sum_{p \in P_w} x_p^* \ln\left( x_p^* / d_w \right) \right) \le$$

$$\le \left\{ \gamma\psi\left( t^N/\gamma \right) + \sum_{e \in E} \sigma_e^*\left( t_e^N \right) \right\} + \left\{ \sum_{e \in E} \sigma_e\left( \bar{f}_e^N \right) + \gamma \sum_{w \in OD} \sum_{p \in P_w} \bar{x}_p^N \ln\left( \bar{x}_p^N / d_w \right) \right\} \le \frac{\tilde{R}^2}{A_N} + \frac{\varepsilon}{2}, \quad (15)$$

где $\left( f^*, x^* \right)$ – решение задачи (1).

На формулу (14) можно смотреть как на критерий останова УМПТ. А именно, ждем когда (вычислимый с помощью быстрого автоматического дифференцирования п. 6) зазор двойственности станет меньше $\varepsilon$. Формула (10) (с заменой $R$ на $\tilde{R}$) содержит гарантированную оценку, на число итераций УМПТ, после которого метод гарантированно остановится по критерию (14).

Для модели (стохастической) стабильной динамики ($\mu \to 0+$) необходимо немного по-другому провести рассуждения

$$\gamma\psi\left( t^N/\gamma \right) + \left\langle \bar{f}, t^N - \bar{t} \right\rangle \le$$

$$\le \min_{t \ge \bar{t}} \left\{ \frac{1}{A_N} \left[ \sum_{k=0}^{N} \alpha_k \cdot \left( \gamma\psi\left( y^k/\gamma \right) + \left\langle \nabla\gamma\psi\left( y^k/\gamma \right), t - y^k \right\rangle \right) \right] + \left\langle \bar{f}, t - \bar{t} \right\rangle + \frac{1}{2A_N} \left\| t - \bar{t} \right\|_2^2 \right\} + \frac{\varepsilon}{2} \le$$

$$\le \min_{t \ge \bar{t},\, \left\| t - \bar{t} \right\|_2^2 \le 2R^2} \left\{ \frac{1}{A_N} \left[ \sum_{k=0}^{N} \alpha_k \cdot \left( \gamma\psi\left( y^k/\gamma \right) + \left\langle \nabla\gamma\psi\left( y^k/\gamma \right), t - y^k \right\rangle \right) \right] + \left\langle \bar{f}, t - \bar{t} \right\rangle \right\} + \frac{R^2}{A_N} + \frac{\varepsilon}{2},$$

где $R^2 = \frac{1}{2} \left\| t_* - y^0 \right\|_2^2 = \frac{1}{2} \left\| t_* - \bar{t} \right\|_2^2$. Следовательно

$$\gamma\psi\left( t^N/\gamma \right) + \left\langle \bar{f}, t^N - \bar{t} \right\rangle -$$

$$-\min_{t\ge\bar t,\,\|t-\bar t\|_2^2\le 10R^2}\left\{\frac{1}{A_N}\left[\sum_{k=0}^{N}\alpha_k\cdot\left(\gamma\psi\left(y^k/\gamma\right)+\left\langle\nabla\gamma\psi\left(y^k/\gamma\right),t-y^k\right\rangle\right)\right]+\left\langle\bar f,t-\bar t\right\rangle\right\}\le\frac{5R^2}{A_N}+\frac{\varepsilon}{2}.$$

Поскольку

$$-\min_{t\ge\bar t,\,\|t-\bar t\|_2^2\le 10R^2}\left\{\frac{1}{A_N}\left[\sum_{k=0}^{N}\alpha_k\left\langle\nabla\gamma\psi\left(y^k/\gamma\right),t\right\rangle\right]+\left\langle\bar f,t-\bar t\right\rangle\right\}=$$
$$=\max_{t\ge\bar t,\,\|t-\bar t\|_2^2\le 10R^2}\left\{\left\langle\frac{1}{A_N}\sum_{k=0}^{N}\alpha_k f^k,t\right\rangle-\left\langle\bar f,t-\bar t\right\rangle\right\}=$$
$$=\max_{t\ge\bar t,\,\|t-\bar t\|_2^2\le 10R^2}\left\{\left\langle\bar f^N-\bar f,t-\bar t\right\rangle\right\}+\left\langle\bar f^N,\bar t\right\rangle\ge$$
$$\ge\left\langle\bar f^N,\bar t\right\rangle+3R\left\|\left(\bar f^N-\bar f\right)_+\right\|_2,$$

то

$$\gamma\psi\left(t^N/\gamma\right)+\left\langle\bar f,t^N-\bar t\right\rangle+\left\langle\bar f^N,\bar t\right\rangle+\gamma\sum_{w\in OD}\sum_{p\in P_w}\bar x_p^N\ln\left(\bar x_p^N/d_w\right)+3R\left\|\left(\bar f^N-\bar f\right)_+\right\|_2\le$$
$$\gamma\psi\left(t^N/\gamma\right)+\left\langle\bar f,t^N-\bar t\right\rangle+\left\langle\bar f^N,\bar t\right\rangle+\gamma\frac{1}{A_N}\sum_{k=0}^{N}\alpha_k\sum_{w\in OD}\sum_{p\in P_w}x_p^k\ln\left(x_p^k/d_w\right)+3R\left\|\left(\bar f^N-\bar f\right)_+\right\|_2\le$$
$$\le\frac{5R^2}{A_N}+\frac{\varepsilon}{2}. \qquad (16)$$

Повторяя рассуждения п. 3 [33] и п. 6.11 [36] из (16), получим

$$\left\langle\bar f^N,\bar t\right\rangle+\gamma\sum_{w\in OD}\sum_{p\in P_w}\bar x_p^N\ln\left(\bar x_p^N/d_w\right)-\left(\left\langle f^*,\bar t\right\rangle+\gamma\sum_{w\in OD}\sum_{p\in P_w}x_p^*\ln\left(x_p^*/d_w\right)\right)\le$$
$$\le\gamma\psi\left(t^N/\gamma\right)+\left\langle\bar f,t^N-\bar t\right\rangle+\left\langle\bar f^N,\bar t\right\rangle+\gamma\sum_{w\in OD}\sum_{p\in P_w}\bar x_p^N\ln\left(\bar x_p^N/d_w\right)\le\frac{5R^2}{A_N}+\frac{\varepsilon}{2}, \qquad (17)$$

$$\left\|\left(\bar f^N-\bar f\right)_+\right\|_2\le\frac{5R}{A_N}+\frac{\varepsilon}{2R}, \qquad (18)$$

где $\left(f^*,x^*\right)$ – решение задачи (2).

На формулу (16) можно смотреть как на критерий останова УМПТ. А именно, ждем когда (вычислимая) левая часть (16) станет меньше $\varepsilon$. Формула (10) (с заменой $R$ на $\sqrt{5}R$) содержит гарантированную оценку, на число итераций УМПТ, после которого метод гарантированно остановится по критерию (16).

Далее для краткости будем единым образом обозначать через $\bar R$ либо $\tilde R$, либо $\sqrt{5}R$. Расшифровывать $\bar R$ нужно будет в зависимости от контекста. Если $\mu\to 0+$, то $\bar R=\sqrt{5}R$, иначе $\bar R=\tilde R$.

Все приведенные в этом пункте формулы выдерживают предельные переходы $\gamma\to 0+$, $\mu\to 0+$ (см. п. 2).

## 6. Вычисление (суб-)градиентов в задаче поиска равновесий в транспортных сетях

Приведем, следуя [5, 27], сглаженный идемпотентный аналог метода Форда–Беллмана [29, 41], позволяющий эффективно рассчитывать значение характеристической функции $\gamma\psi\left(t/\gamma\right)$. Для этого предположим, что любые движения по ребрам графа с учетом их ориентации являются допустимыми, т.е. множество путей, соединяющих заданные две вершины (источник и сток), – это множество всевозможных способов добраться из источника в сток по ребрам имеющегося графа с учетом их ориентации. Этого всегда можно добиться раздутием исходного графа в несколько раз за счет введения дополнительных вершин и ребер. Такое раздутие заведомо можно сделать за $\mathrm{O}(n)$. Отметим, что при этом в качестве путей будут присутствовать, в том числе, и самопересекающиеся маршруты. Однако можно показать, что вклад таких "не физических" путей в итоговую равновесную конфигурацию будет пренебрежимо мал.

Будем считать, как и прежде, что число ребер в любом пути не больше $H = \mathrm{O}\left(\sqrt{n}\right)$. Введем классы путей: $P_{ij}^l$ – множество всех путей из $i$ в $j$, состоящих ровно из $l$ ребер, $\tilde{P}_{ij}^l$ – множество всех путей из $i$ в $j$, состоящих из не более чем $l$ ребер. Зафиксируем источник (вершину) $i \in V$, и введем следующие функции для $j \in V$, $l = 1,\dots,H$ :

$$\begin{cases} a_{ij}^l(t) = \gamma\psi_{P_{ij}^l}(t/\gamma) = \gamma \ln\left(\sum\limits_{p\in P_{ij}^l} \exp\left(-\sum\limits_{e\in E}\delta_{ep}t_e \Big/ \gamma\right)\right) \\ b_{ij}^l(t) = \gamma\psi_{\tilde{P}_{ij}^l}(t/\gamma) = \gamma \ln\left(\sum\limits_{p\in \tilde{P}_{ij}^l} \exp\left(-\sum\limits_{e\in E}\delta_{ep}t_e \Big/ \gamma\right)\right) \end{cases}.$$

Некоторые из этих функций могут быть равны $-\infty$. Это означает, что соответствующее множество маршрутов – пустое. Данные функции можно вычислять рекурсивным образом:

$$a_{ij}^1(t) = b_{ij}^1(t) = \begin{cases} -t_e,\ e = (i \to j) \in E \\ -\infty,\ e = (i \to j) \notin E \end{cases},$$

$$\begin{cases} a_{ij}^{l+1}(t) = \gamma \ln\left(\sum\limits_{k:\, e=(k\to j)\in E} \exp\left(\left(a_{ik}^l(t) - t_e\right)/\gamma\right)\right) \\ b_{ij}^{l+1}(t) = \gamma\ln\left(\exp\left(b_{ij}^l(t)/\gamma\right) + \exp\left(a_{ij}^{l+1}(t)/\gamma\right)\right) \end{cases},\ j \in V,\ l = 1,\dots,H-1. \quad (19)$$

На каждом шаге $l$ необходимо сделать $\mathrm{O}(n)$ арифметических операций. Следовательно, для вычисления $\gamma\psi(t/\gamma)$ необходимо сделать $\mathrm{O}(SHn)$ арифметических операций, где $S = |O|$ – число источников, как правило, можно считать $S \ll n$. Причем вычисление функции $\gamma\psi(t/\gamma)$ (и ее градиента) может быть распараллелено на $S$ процессорах. При $\gamma \to 0+$ процедура вырождается в известный метод Форда–Беллмана [29, 41] (динамическое программирование). В процедуре Форда–Беллмана требуется посчитать $H$-степень матрицы $A = \left\|a_{ij}\right\|_{i,j\in V}$,

$$a_{ij} = t_e,\ e = (i \to j) \in E;$$
$$a_{ij} = \infty,\ e = (i \to j) \notin E,$$

в идемпотентной математике (вместо обычного поля используется тропическое полуполе [42] со следующими операциями: сложение $a \oplus b = \min\{a,b\}$, произведение $a \otimes b = a + b$). Учитывая, что

$$a \oplus b = \min\{a,b\} = -\lim_{\gamma\to 0+} \gamma\ln\left(\exp(-a/\gamma) + \exp(-b/\gamma)\right),$$

$$a \otimes b = a + b = -\lim_{\gamma \to 0+} \gamma \ln\left(\exp(-a/\gamma) \cdot \exp(-b/\gamma)\right),$$

можно посчитать обычную (над обычным полем) $H$ -степень матрицы $A^\gamma = \left\| a_{ij}^\gamma \right\|_{i,j \in V}$,

$$a_{ij}^\gamma = e^{-t_e/\gamma},\ e = (i \to j) \in E;$$

$$a_{ij}^\gamma = 0,\ e = (i \to j) \notin E,$$

и применить поэлементно к полученной матрице $-\gamma \ln(\,\cdot\,)$. В пределе $\gamma \to 0+$ получим метод Форда–Беллмана. Однако если не делать предельный переход, то получается нужный нам сглаженный вариант этого алгоритма с такой же временной сложностью. Некоторый аналог этого сглаженного варианта, по сути, и был описан выше.

Используя быстрое автоматическое дифференцирования [28], опишем способ вычисления градиента функции $\gamma\psi^i(t/\gamma) = \sum_{j:\, w=(i,j) \in OD} d_w \gamma \psi_w(t/\gamma)$:

$$\frac{\partial \psi^i}{\partial b_{ij}^H} = d_{ij},\ \frac{\partial \psi^i}{\partial a_{ij}^H} = 0,$$

$$\begin{cases} \dfrac{\partial \psi^i}{\partial b_{ij}^l} = \dfrac{\partial \psi^i}{\partial b_{ij}^{l+1}} \dfrac{\partial b_{ij}^{l+1}}{\partial b_{ij}^l} \\ \dfrac{\partial \psi^i}{\partial a_{ij}^l} = \dfrac{\partial \psi^i}{\partial b_{ij}^l} \dfrac{\partial b_{ij}^l}{\partial a_{ij}^l} + \sum\limits_{k:\, e=(j \to k) \in E} \dfrac{\partial \psi^i}{\partial a_{ik}^{l+1}} \dfrac{\partial a_{ik}^{l+1}}{\partial a_{ij}^l} \end{cases},\ j \in V,\ l = H-1, \ldots, 1. \tag{20}$$

$$\frac{\partial \psi^i}{\partial t_e} = \sum_{l=0}^{H-1} \sum_{j \in V} \frac{\partial \psi^i}{\partial a_{ij}^{l+1}} \frac{\partial a_{ij}^{l+1}}{\partial t_e} = \sum_{l=0}^{H-1} \frac{\partial \psi^i}{\partial a_{ij'}^{l+1}} \frac{\partial a_{ij'}^{l+1}}{\partial t_e},\ e = (k, j'). \tag{21}$$

Частные производные

$$\frac{\partial b_{ij}^{l+1}}{\partial b_{ij}^l},\ \frac{\partial b_{ij}^l}{\partial a_{ij}^l},\ \frac{\partial a_{ik}^{l+1}}{\partial a_{ik}^l},\ \frac{\partial a_{ij}^{l+1}}{\partial t_e}$$

могут быть явно вычислены из системы (19) за $O(1)$ каждая. Поясним примером. Пусть, например, для некоторых $l$ и $j$ имеет место

$$a_{ij}^{l+1}(t) = \gamma \ln\left(\exp\left(\left(a_{ik_1}^l(t) - t_{e_1}\right)/\gamma\right) + \exp\left(\left(a_{ik_2}^l(t) - t_{e_2}\right)/\gamma\right)\right).$$

Тогда[7]

$$\frac{\partial a_{ij}^{l+1}}{\partial a_{ik_1}^l} = \frac{\exp\left(\left(a_{ik_1}^l(t) - t_{e_1}\right)/\gamma\right)}{\exp\left(\left(a_{ik_1}^l(t) - t_{e_1}\right)/\gamma\right) + \exp\left(\left(a_{ik_2}^l(t) - t_{e_2}\right)/\gamma\right)},$$

---

[7] Заметим, что возникающие при расчете градиента отношения экспонент стоит сразу же приводить (для большей вычислительной устойчивости) к дробям с числителем равным 1

$$\frac{e^a}{e^a + e^b + \ldots} = \frac{1}{1 + e^{b-a} + \ldots}.$$

Аналогично с выражениями вида

$$\ln\left(\exp(a) + \exp(b)\right) = c + \ln\left(\exp(a-c) + \exp(b-c)\right),\ c = \max\{a, b\}.$$

$$\frac{\partial a_{ij}^{l+1}}{\partial t_{e_1}} = -\frac{\exp\left(\left(a_{ik_1}^{l}(t)-t_{e_1}\right)/\gamma\right)}{\exp\left(\left(a_{ik_1}^{l}(t)-t_{e_1}\right)/\gamma\right)+\exp\left(\left(a_{ik_2}^{l}(t)-t_{e_2}\right)/\gamma\right)}.$$

Остальные частные производные являются неизвестными. Несложно понять, что система (20) может быть последовательно разрешена за $\mathrm{O}(Hn)$. Фактически тут происходит процесс прохождения графа вычислений (19) в обратном направлении, с той лишь разницей, что обратный проход получается приблизительно в 4 раза дороже. К сожалению, в отличие от прямого процесса (19), теперь, чтобы вычислить градиент, необходимо хранить в памяти промежуточные вычисления, см. формулу (21). То есть для вычисления градиента $\gamma\psi^i(t/\gamma)$ требуется время $\mathrm{O}(Hn)$ и память $\mathrm{O}(Hn)$. Таким образом, для вычисления градиента $\gamma\psi(t/\gamma)=\sum_{i\in O}\gamma\psi^i(t/\gamma)$ требуется время $\mathrm{O}(SHn)$ и память $\mathrm{O}(SHn)$. Также как и при вычислении значения функции $\gamma\psi(t/\gamma)$, описанная выше схема вычисления градиента может быть распараллелена на $S=|O|$ процессорах.

Заметим, что расчет $\nabla\gamma\psi(t/\gamma)$, реализованный на языке программирования MatLab 8 на стандартном современном ноутбуке (с тактовой частотой 1.9 ГГц), для манхетенского графа с числом ребер $n \simeq 10^3$, занимал порядка одной минуты, в то время как использование прямого дифференцирования занимало порядка 100 минут.

В связи с написанным выше, важно отметить, в пределе $\gamma \to 0+$ оценка $\mathrm{O}(SHn)$ сложности вычисления градиента $\gamma\psi(t/\gamma)$ может быть улучшена до оценки $\mathrm{O}(Sn\ln n)$ вычисления субградиента (5), (6) [6, 29]. Действительно, для каждого из $S$ источников можно построить (например, алгоритмом Дейкстры [29]) соответствующее дерево кратчайших путей. Исходя из принципа динамического программирования “часть кратчайшего пути сама будет кратчайшим путем” несложно понять, что получится именно дерево, с корнем в рассматриваемом источнике. Это можно сделать для одного источника за $\mathrm{O}(n\ln n)$ [6, 12, 29]. Однако, главное, правильно взвешивать ребра (их $n$) такого дерева, чтобы за один проход этого дерева можно было восстановить вклад (по всем ребрам) соответствующего источника в субградиент. Ребро должно иметь вес равный сумме всех проходящих через него корреспонденций с заданным источником (корнем дерева). Имея значения соответствующих корреспонденций (их не больше $\mathrm{O}(n)$) за один обратный проход (то есть с листьев к корню) такого дерева можно осуществить необходимое взвешивание (с затратами не более $\mathrm{O}(n)$). Делается это по правилу: вес ребра равен сумме корреспонденции (возможно, равной нулю), в соответствующую вершину, в которую ребро входит и сумме весов всех ребер (если таковые имеются), выходящих из упомянутой вершины.

## 7. Сопоставление оценок, численные эксперименты

Используя наработки предыдущих пунктов, в частности, формулы (10)–(12), (15), (17), (18), можно резюмировать полученные в статье результаты в виде следующей теоремы.

**Теорема.** *УМПТ из п. 4 с критериями останова (14), (16) гарантированно остановится, достигнув желаемой точности* $\varepsilon$*, сделав арифметических операций не больше, чем приведено в соответствующем поле таблицы 1.*

Таблица 1

| Время работы | $\gamma > 0$ | $\gamma \to 0+$ |
|---|---|---|
| $\mu \to 0+$, $\mu > 0$ – не важно | $\mathrm{O}\left( SHn \cdot \sqrt{\dfrac{Hd\bar{R}^2}{\gamma\varepsilon}} \right)$ (22) | $\mathrm{O}\left( Sn\ln n \cdot \dfrac{Hd^2\bar{R}^2}{\varepsilon^2} \right)$ (23) |

Обратим внимание, что оценки (22), (23) не зависят от (потенциально экспоненциально большого) числа путей $|P|$. В частности, $|P| \gg 2^{\sqrt{n}}$, для манхетенских сетей. Также обратим внимание, что оценка (22) весьма чувствительна к предельному переходу $\gamma \to 0+$.

Важно также заметить, что обычные (не стохастические) равновесия можно искать с помощью искусственного введения энтропийной регуляризации [9, 17]. При этом, чтобы с точностью $\varepsilon > 0$ по функции решить исходную задачу (1) или (2) с $\gamma = 0$, можно действовать следующим образом. Выбрать

$$\gamma_* = \frac{\varepsilon}{2\sum_{w \in OD} d_w \ln|P_w|}$$

и решать регуляризованную задачу (1) или (2) с $\gamma = \gamma_*$, но с точностью $\varepsilon/2$ [17]. В этом случае оценка (22) из таблицы 1 примет вид

$$\mathrm{O}\left( SHn \cdot \frac{d}{\varepsilon}\sqrt{H\bar{R}^2\omega} \right), \tag{24}$$

где $\omega = \max_{w \in OD} \ln|P_w|$ (в худшем случае можно ожидать $\omega \sim \sqrt{n}$). Сопоставляя оценки (23), (24) можно прийти к выводу, что введенная регуляризация (приводящая к сглаживанию двойственной задачи (3) [38]) оправдана. Однако не стоит забывать, что все приведенные оценки – это верхние оценки. Численные эксперименты показали, что эти оценки не всегда точны, особенно в части формулы (23).

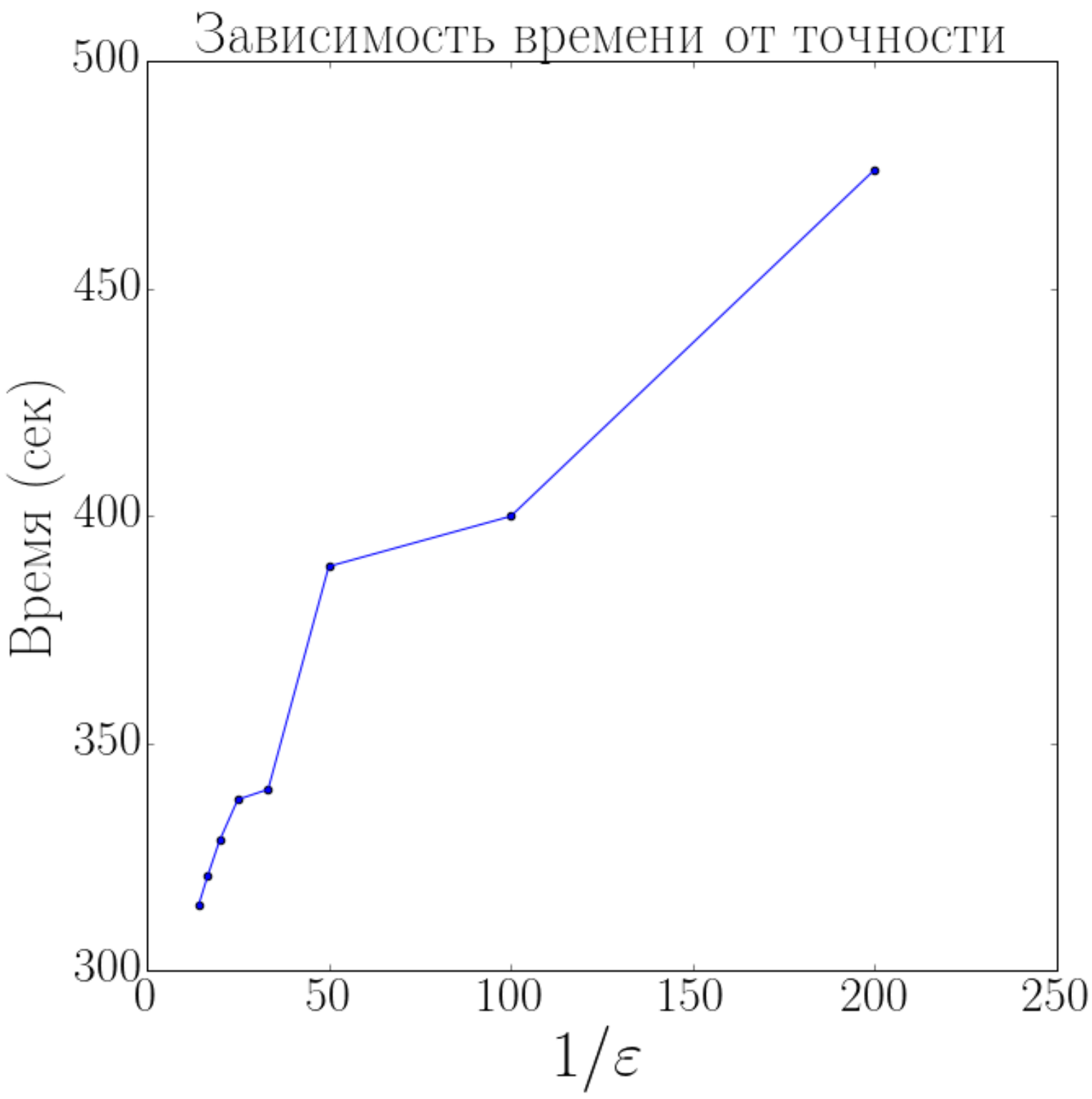

Рис. 1 Модель Бэкмана ($\mu = 0.25$), $\gamma \to 0+$ ($\tilde{\varepsilon}$ обозначено здесь через $\varepsilon$)

Прежде всего, заметим, что оценки общего числа арифметических операций из таблицы 1 можно понимать как [стоимость итерации]×[число итераций]. А число итераций $N$, с некоторой натяжкой, можно просто считать пропорциональным относительной точности $\tilde{\varepsilon}$ в некоторой степени[8] $N \sim \tilde{\varepsilon}^{-\beta}$. Именно такая (относительная) точность, как правило, интересна на практике. Из теоретических верхних оценок (22)–(24) напрашивается вывод, что при $\gamma \gg \gamma_*$ имеем $N \sim \tilde{\varepsilon}^{-1/2}$, при $\gamma \simeq \gamma_*$ имеем $N \sim \tilde{\varepsilon}^{-1}$, а при $\gamma \to 0+$ имеем $N \sim \tilde{\varepsilon}^{-2}$. На самом деле, конечно, опущенные числовые множители тут также играют важную роль. Однако еще более важно то, что численные эксперименты [43] не подтвердили оценку $N \sim \tilde{\varepsilon}^{-2}$ при $\gamma \to 0+$. Более того, наблюдалась совсем другая зависимость $N \sim \tilde{\varepsilon}^{-1}$ [43], точнее говоря, наблюдалась зависимость $N \sim C_1 + C_2\tilde{\varepsilon}^{-1}$ (см. рис. 1, 2). Учитывая, что стоимость итерации заметно меньше в случае $\gamma \to 0+$ ($SHn \to Sn\ln n$, см. п. 6), то вывод об оправданности регуляризации приходится поставить под сомнение. Для графа города Анахайм (Anaheim) [44] с $n \simeq 10^3$, $S \sim 40$ и $|OD| \sim 1.5 \cdot 10^3$ УМПТ, реализованный на стандартном современном ноутбуке (с тактовой частотой 1.9 ГГц) на языке программирования Python 2.4 [43], находил с относительной точностью $\tilde{\varepsilon} \simeq 0.01$ равновесие в модели Бэкмана ($\mu > 0, \gamma \to 0+$) приблизительно за 7 минут, в модели стабильной динамики ($\mu > 0, \gamma \to 0+$) за 10 минут (см. рис. 1, 2). Для соответствующей регуляризованной модели – заметно дольше [45].

---

[8] $\tilde{\varepsilon} = 1\% = 0.01$ – означает, что начальную невязку по функции или зазору двойственности необходимо уменьшить в 100 раз.

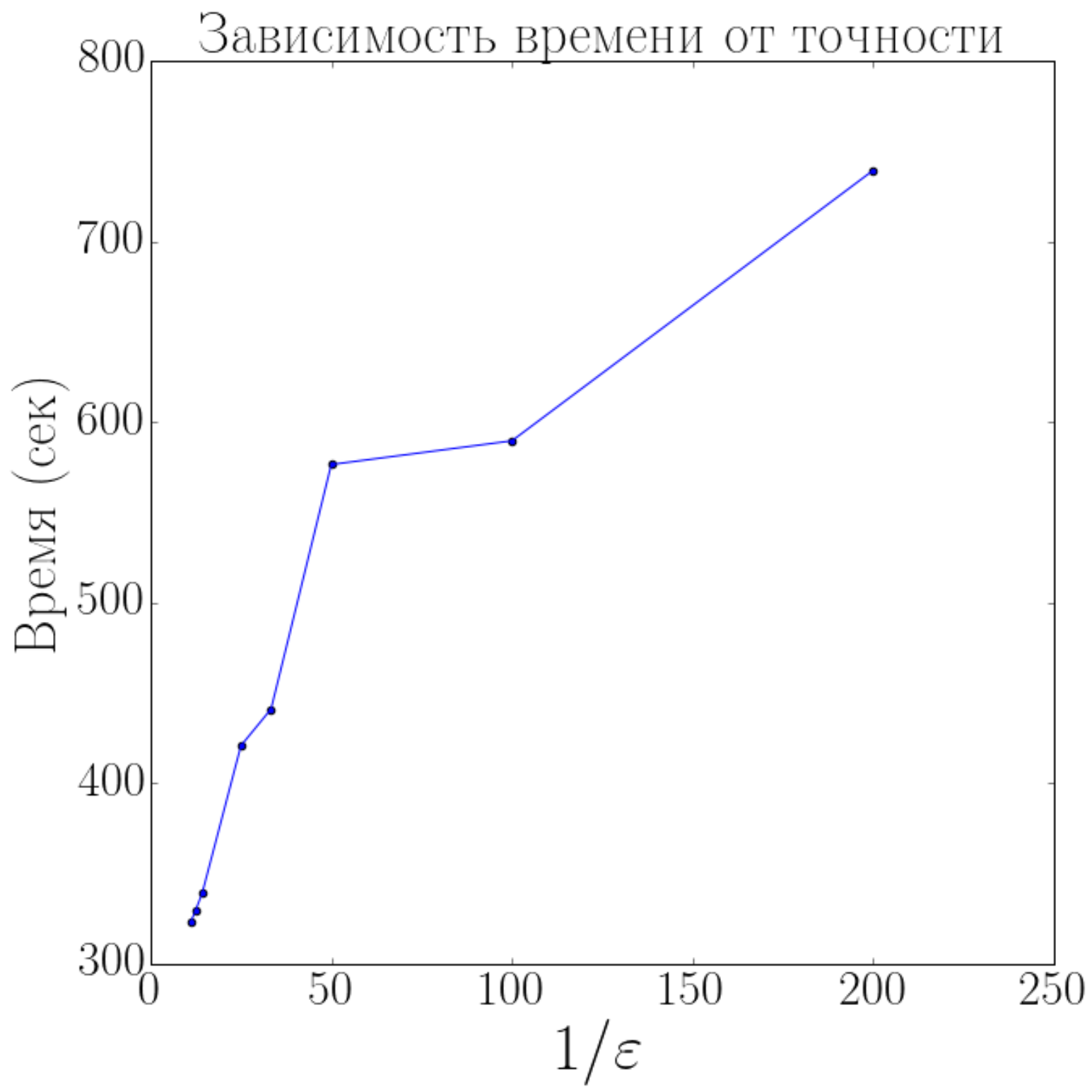

Рис. 2 Модель стабильной динамики ( $\mu \to 0+$ ), $\gamma \to 0+$ ( $\tilde{\varepsilon}$ обозначено здесь через $\varepsilon$ )

Отметим, что написанный код не был оптимизирован. Мы полагаем, что при правильной реализации и с использованием языка более низкого уровня (например, С++) можно ожидать выигрыша как минимум на порядок (скорее всего, на два). В таком случае можно рассчитывать, что предлагаемый в данной статье метод приблизится по производительности к методу условного градиента [3, 6] (Франк–Вульфа), который считается сейчас наиболее эффективным методом поиска равновесного распределения потоков по путям (ребрам) в модели Бэкмана ( $\mu = 0.25$, $\gamma \to 0+$ ). В лучшей (из известных нам) реализаций этого метода, сделанной А.С. Аникиным и А.Ю. Горновым на С++ [46], для графа, с аналогичными параметрами и с аналогичными требованиями к точности, метод условного градиента сходился за пару секунд. Для модели стабильной динамики имеющиеся сейчас реализации (при $\gamma \to 0+$ ) существующих методов [6, 12] заметно проигрывают, изложенным в данной работе.

Описанный в статье подход можно переписать, беря в качестве базисного метода УМПТ из статьи [17], в варианте, работающем с сильно выпуклыми постановками задач. Поскольку двойственная задача (3) не является сильно выпуклой (во всяком случае, на данный момент умеют устанавливать только выпуклость этой задачи), то необходимо регуляризовывать двойственную задачу. Можно показать, по-сути, рассуждая подобно статье [39] (см. также концовку статьи [33]), что при правильной регуляризации существенно упрощаются формулы восстановления решения прямой задачи.

Описанный в статье формализм может быть перенесен и на задачи поиска равновесий в многостадийных транспортных моделях [1, 7, 8, 10, 11, 13, 30]. Для этого стоит использовать подход работ [11, 13]. Этому со временем планируется посвятить отдельную публикацию.

Авторы выражают благодарность Ю.Е. Нестерову за ряд ценных замечаний.

Исследование в пп. 4, 5, 6 выполнено ИППИ РАН за счет гранта Российского научного фонда (проект № 14-50-00150). Исследование в других пунктах выполнено ИППИ РАН при поддержке гранта РФФИ 15-31-70001-мол_а_мос. Исследование П.Е. Двуреченского и М.Б. Кубентаевой в п. 4 также было поддержано грантом Президента РФ МК-1806.2017.9.

## Литература